\documentclass{amsart}
\usepackage{amsmath}
\usepackage{amssymb}

\newtheorem{theorem}{Theorem}
\newtheorem{lemma}[theorem]{Lemma}

\newtheorem{definition}[theorem]{Definition}

\def\endproof{\hfill\qedbox\smallbreak\noindent}
\def\qedbox{\hbox{$\rlap{$\sqcap$}\sqcup$}}
\def\MM{\mathfrak{M}}\def\NN{\mathfrak{N}}
\usepackage{color}

\begin{document}
  \title[Simple Jacobi--Ricci commuting tensors]{{The classification of simple Jacobi--Ricci commuting algebraic
curvature tensors}}
\author{P. Gilkey}
  \begin{address}{Mathematics Department, University of Oregon,
 Eugene Oregon 97403 USA}\end{address}
    \begin{email}{gilkey@uoregon.edu}\end{email}
\author{S. Nik\v cevi\'c}
 \begin{address}{Mathematical Institute, SANU,
 Knez Mihailova 35, p.p. 367,
 11001 Belgrade,
 Serbia.}\end{address}
    \begin{email}{stanan@mi.sanu.ac.yu}\end{email}
\begin{abstract}
We classify algebraic curvature tensors such that the Ricci operator $\rho$ is simple
(i.e. $\rho$ is complex diagonalizable and either $\operatorname{Spec}\{\rho\}=\{a\}$ or
$\operatorname{Spec}\{\rho\}=\{a_1\pm a_2\sqrt{-1}\}$) and which are Jacobi--Ricci commuting (i.e. 
$\rho\mathcal{J}(v)=\mathcal{J}(v)\rho$ for all $v$).
\end{abstract}
 \maketitle
\centerline{This is dedicated to Professor Oldrich Kowalski}\bigbreak
The study of curvature is fundamental in differential geometry. It is often convenient to work first in an abstract
algebraic context and then subsequently to pass to the geometrical setting. We say that
$\MM:=(V,\langle\cdot,\cdot\rangle,A)$ is a {\it model} if $\langle\cdot,\cdot\rangle$ is a non-degenerate
inner-product of signature $(p,q)$ on a real vector space $V$ of dimension $m=p+q$ and if $A\in\otimes^4V^*$ is an
{\it algebraic curvature tensor}, i.e. a $4$-tensor which has the symmetries of the Riemann curvature tensor:
\begin{equation}\label{eqn-1}
\begin{array}{l}
A(v_1,v_2,v_3,v_4)=-A(v_2,v_1,v_3,v_4)=A(v_3,v_4,v_1,v_2),\\
A(v_1,v_2,v_3,v_4)+A(v_2,v_3,v_1,v_4)+A(v_3,v_1,v_2,v_4)=0\,.
\end{array}\end{equation}
If $P$ is a point of a pseudo-Riemannian manifold $\mathcal{M}=(M,g)$, then the associated model is defined by
setting $\MM(\mathcal{M},P):=(T_PM,g_P,R_P)$ where $R_P$ is the curvature tensor of the Levi--Civita connection;
every model is geometrically realizable in this fashion. Consequently the study of algebraic curvature tensors plays a
central role in many geometric investigations.

If $\MM$ is a model, then {\it Jacobi operator}
$\mathcal{J}$, the skew-symmetric curvature operator $\mathcal{R}$, and the {\it Ricci operator} $\rho$
are defined by the identities:
\begin{eqnarray*}
&&
\langle\mathcal{J}(x)y,z\rangle=A(y,x,x,z),\qquad
  \langle\mathcal{R}(x,y)z,w\rangle=A(x,y,z,w),\\
&&\langle\rho x,y\rangle
  =\operatorname{Tr}\left\{z\rightarrow\textstyle\frac12\mathcal{R}(z,x)y+\frac12\mathcal{R}(z,y)x\right\}\,.
\end{eqnarray*}
One says $\MM$ is {\it Einstein}  if $\rho=a\operatorname{id}$; $a$ is called the {\it Einstein constant}.

The study of commutativity properties of natural operators defined by the curvature tensor was initiated by
Stanilov \cite{S1,S2} and has proved to be a very fruitful one; we refer to \cite{B07} for a survey of the field
and for a more complete bibliography than is possible to present here. 

We begin with the following fundamental result which is established in \cite{GN07} and which examines when the
Jacobi operator or the skew-symmetric curvature operator commutes with the Ricci operator:

\begin{lemma}\label{lem-1}
The following conditions are equivalent for a model $\MM$:
\begin{enumerate}
\item $\mathcal{J}(v)\rho=\rho\mathcal{J}(v)$ for all $v\in V$.
\item $\mathcal{R}(v_1,v_2)\rho=\rho\mathcal{R}(v_1,v_2)$ for all $v_1,v_2\in V$.
\item $A(\rho v_1,v_2,v_3,v_4)=A(v_1,\rho v_2,v_3,v_4)=A(v_1,v_2,\rho
v_3,v_4)$\newline$\phantom{.......}=A(v_1,v_2,v_3,\rho v_4)$ for all
$v_1,v_2,v_3,v_4\in V$.
\end{enumerate}
\end{lemma}

One says that a model $\MM$ is {\it decomposable} if there is an orthogonal direct sum decomposition
$V=V_1\oplus V_2$ inducing a splitting $A=A_1\oplus A_2$; $\MM$ is {\it indecomposable} if it is not decomposable. Let
$\operatorname{Spec}(\rho)\subset\mathbb{C}$ be the spectrum of the Ricci operator.  One has \cite{GPV07}:

\begin{lemma}\label{lem-2}
If $\MM$ is an indecomposable Jacobi--Ricci commuting model, then either $\operatorname{Spec}(\rho)=\{a_1\}$ 
or $\operatorname{Spec}(\rho)=\{a_1\pm a_2\sqrt{-1}\}$ where $a_2>0$.
\end{lemma}

Although $\rho$ is self-adjoint, $\rho$ need not be diagonalizable in the higher signature setting and in fact the
Jordan normal form of $\rho$ can be quite complicated. To simplify the discussion, we shall suppose $\rho$ complex
diagonalizable henceforth. Motivated by Lemmas \ref{lem-1} and \ref{lem-2}, we make the following:

\begin{definition}\label{defn-3} 
If $\MM$ is a model which has any of the (3) equivalent properties listed in
Lemma
\ref{lem-1}, then $\MM$ is said to be a {\it Jacobi--Ricci commuting model}. If in addition, the Ricci operator $\rho$
is complex diagonalizable and if either $\operatorname{Spec}(\rho)=\{a_1\}$ or $\operatorname{Spec}(\rho)=\{a_1\pm
a_2\sqrt{-1}\}$ for
$a_2>0$, then $\MM$ is said to be a {\it simple Jacobi--Ricci commuting model}.  If $\mathcal{M}$ is a
pseudo-Riemannian manifold, then $\mathcal{M}$ is said to be {\it simple Jacobi--Ricci commuting}  if
$\MM(\mathcal{M},P)$ is a simple Jacob--Ricci commuting model for all points $P$ of $M$. 

\end{definition}

The following is immediate from the definitions we have given.
\begin{lemma}\label{lem-4}
Let $\MM$ be a model.
\begin{enumerate}
\item $\MM$ is Einstein if and only if $\MM$ is a simple Jacobi--Ricci commuting model with
$\operatorname{Spec}(\rho)=\{a_1\}$.
\item Let $\MM$ be a simple Jacobi--Ricci commuting model which is not Einstein. Set
$J=J_\MM:=a_2^{-1}\{\rho-a_1\operatorname{id}\}$. Then $J$ is a self-adjoint complex structure on $V$,
$A(Jx,y,z,w)=A(x,Jy,z,w)=A(x,y,Jz,w)=A(x,y,z,Jw)$, and $\rho=a_1+a_2J$.
\end{enumerate}\end{lemma}

In view of Lemma \ref{lem-4} (1), we shall assume $\MM$ is not Einstein henceforth. The following ansatz for
constructing simple Jacobi--Ricci commuting models which are not Einstein will be crucial:

\begin{definition}\label{defn-5}
\rm Let $\NN:=(V_0,g,A_1,A_2)$ where $g$ is a positive definite inner
product on a finite dimensional real vector space $V_0$ and where $A_1$ and $A_2$ are Einstein algebraic curvature
tensors with Einstein constants, respectively, $a_1$ and $a_2>0$. Extend
$g$, $A_1$, and $A_2$ to be complex linear on the complexification
$V_{\mathbb{C}}:=V_0\otimes_{\mathbb{R}}\mathbb{C}$. Let
$V:=V_0\oplus V_0\sqrt{-1}$ be the underlying real vector space of $V_{\mathbb{C}}$, let
$\langle\cdot,\cdot\rangle:=\operatorname{Re}g(\cdot,\cdot)$, and let
$A:=\operatorname{Re}\{A_1+\sqrt{-1}A_2\}$ define $\MM(\NN):=(V,\langle\cdot,\cdot\rangle,A)$.
\end{definition}

The following classification result is the fundamental result of this paper:

\begin{theorem}\label{thm-6}
Adopt the notation established above:
\begin{enumerate}
\item $\MM(\NN)$ is a simple Jacobi--Ricci commuting model which is not Einstein, which has
$\operatorname{Spec}\{\rho\}=\{2a_1\pm2a_2\sqrt{-1}\}$, and which has that $J_{\MM(\NN)}$ is multiplication by
$\sqrt{-1}$.
\item Let $\MM$ be a simple Jacobi--Ricci commuting model which is not Einstein with
$\operatorname{Spec}(\rho)=\{2a_1\pm2a_2\sqrt{-1}\}$ for $a_2>0$. Then there exists
$\NN$ so
$\MM$ is isomorphic to
$\MM(\NN)$.
\item If $\NN=(V_0,g,A_1,A_2)$ and $\tilde\NN=(\tilde V_0,\tilde g,\tilde A_1,\tilde A_2)$, then
$\MM(\NN)$ is isomorphic to
$\MM(\tilde\NN)$ if and only if there is an isomorphism $\theta:V_0\rightarrow\tilde V_0$ and a skew-adjoint linear
transformation $T$ of $(V_0,g)$ with $|T|<1$ so that:
\begin{equation}\label{eqn-2}
\begin{array}{l}
\tilde g(\theta v,\theta w)=g(v,w)-g(Tv,Tw),\\
\tilde A_1(\theta v,\theta w,\theta x,\theta y)=
A_1(v,w,x,y)-A_1(Tv,Tw,x,y)\\
\quad-A_1(Tv,w,Tx,y)-A_1(Tv,w,x,Ty)-A_1(v,Tw,Tx,y)\\
\quad-A_1(v,Tw,x,Ty)-A_1(v,w,Tx,Ty)+A_1(Tv,Tw,Tx,Ty)\\
\quad-A_2(Tv,w,x,y)-A_2(v,Tw,x,y)-A_2(v,w,Tx,y)\\
\quad-A_2(v,w,x,Ty)+A_2(Tv,Tw,Tx,y)+A_2(Tv,Tw,x,Ty)\\
\quad+A_2(Tv,w,Tx,Ty)+A_2(v,Tw,Tx,Ty)\\
\tilde A_2(\theta v,\theta w,\theta x,\theta y)=A_2(v,w,x,y)-A_2(Tv,Tw,x,y)\\
\quad-A_2(Tv,w,Tx,y)-A_2(Tv,w,x,Ty)-A_2(v,Tw,Tx,y)\\
\quad-A_2(v,Tw,x,Ty)-A_2(v,w,Tx,Ty)+A_2(Tv,Tw,Tx,Ty)\\
\quad+A_1(Tv,w,x,y)+A_1(v,Tw,x,y)+A_1(v,w,Tx,y)\\
\quad+A_1(v,w,x,Ty)-A_1(Tv,Tw,Tx,y)-A_1(Tv,Tw,x,Ty)\\
\quad-A_1(Tv,w,Tx,Ty)-A_1(v,Tw,Tx,Ty)\end{array}\end{equation}
\end{enumerate}\end{theorem}

Theorem \ref{thm-6} completes the analysis in the algebraic setting. In the geometric setting, by contrast, the
situation is still far from clear. However there is a geometrical example known
\cite{GN07} in signature
$(2,2)$  which may be described as follows; we refer to \cite{GN07} for further details. Let
$(x_1,x_2,x_3,x_4)$ be coordinates on
$\mathbb{R}^4$. Define a metric whose non-zero components are, up to the usual $\mathbb{Z}_2$ symmetries,
given by:
\begin{equation}\label{eqn-3}
\begin{array}{ll}
g(\partial_1,\partial_3)=g(\partial_2,\partial_4)=1,
&g(\partial_3,\partial_4)=s(x_2^2-x_1^2),\\
g(\partial_3,\partial_3)=2sx_1x_2,\quad
&g(\partial_4,\partial_4)=-2sx_1x_2\,.
\end{array}\end{equation}

\begin{lemma}\label{lem-7}
Let $\mathcal{M}$ be as in Equation (\ref{eqn-3}). Then $\mathcal{M}$ is a locally
symmetric simple Jacobi--Ricci commuting manifold with $\operatorname{Spec}(\rho)=\{\pm
2s\sqrt{-1}\}$ of signature $(2,2)$. The Ricci oprator and non-zero curvatures are described by:
$$\begin{array}{llll}
R_{1314}=s,&R_{1323}=-s,&R_{1424}=s,&R_{2324}=-s,\\
\rho\partial_1=-2s\partial_2,&\rho\partial_2=2s\partial_1,&
\rho\partial_3=2s\partial_4,&\rho\partial_4=-2s\partial_3\,.
\end{array}$$
\end{lemma}

The remainder of this note is devoted to the proof of Theorem \ref{thm-6}. 
\medbreak\noindent{\it Proof of Theorem \ref{thm-6} (1):} We generalize the discussion of
\cite{GN07}. Let $\{e_i\}$ be an orthonormal basis for $V_0$. Let $e_i^+:=e_i$
and
$e_i^-:=\sqrt{-1}e_i$ be an orthonormal basis for $V$; the vectors
$e_i^+$ are spacelike and the vectors $e_i^-$ are timelike so $\MM(\NN)$ has neutral signature. Clearly the
symmetries of Equation (\ref{eqn-1}) hold for the complexification of $A_1$ and $A_2$ and, consequently, for
$A_1+\sqrt{-1}A_2$ and $A=\operatorname{Re}(A_1+\sqrt{-1}A_2)$. Thus
$\MM(\NN)$ is a model and the non-zero components of
$A$ relative to this basis are given by:
\begin{equation}\label{eqn-4}
\begin{array}{l}
A(e_i^-,e_j^+,e_k^+,e_l^+)=A(e_i^+,e_j^-,e_k^+,e_l^+)=A(e_i^+,e_j^+,e_k^-,e_l^+)\\
\quad=A(e_i^+,e_j^+,e_k^+,e_l^-)=-A_2(e_i,e_j,e_k,e_l),\vphantom{\vrule height 11pt}\\
A(e_i^+,e_j^-,e_k^-,e_l^-)=A(e_i^-,e_j^+,e_k^-,e_l^-)=A(e_i^-,e_j^-,e_k^+,e_l^-)\vphantom{\vrule height 11pt}\\
\quad=A(e_i^-,e_j^-,e_k^-,e_l^+)=A_2(e_i,e_j,e_k,e_l)\vphantom{\vrule height 11pt},\\
A(e_i^+,e_j^+,e_k^+,e_l^+)=A(e_i^-,e_j^-,e_k^-,e_l^-)=A_1(e_i,e_j,e_k,e_l),\\
A(e_i^+,e_j^+,e_k^-,e_l^-)=A(e_i^+,e_j^-,e_k^+,e_l^-)=A(e_i^-,e_j^+,e_k^+,e_l^-)
   \phantom{\vrule height 11pt}\\
\quad=A(e_i^+,e_j^-,e_k^-,e_l^+)=A(e_i^-,e_j^+,e_k^-,e_l^+)=A(e_i^-,e_j^-,e_k^+,e_l^+)
\vphantom{\vrule height 11pt}\\
\quad=-A_1(e_i,e_j,e_k,e_l)\,.\vphantom{\vrule height 11pt}
\end{array}\end{equation}
Let $\rho:=\rho_A$ and $\rho_i:=\rho_{A_i}$. We sum over $k$ in the following expansions to see:
\medbreak
$\langle\rho e_i^+,e_j^+\rangle=A(e_i^+,e_k^+,e_k^+,e_j^+)-A(e_i^+,e_k^-,e_k^-,e_j^+)
    =2\rho_1(e_i,e_j)=2a_1\delta_{ij}$,
\par$\langle\rho e_i^-,e_j^-\rangle=A(e_i^-,e_k^+,e_k^+,e_j^-)-A(e_i^-,e_k^-,e_k^-,e_j^-)
    =-2\rho_1(e_i,e_j)=-2a_1\delta_{ij}$,
\par$\langle\rho e_i^+,e_j^-\rangle=A(e_i^+,e_k^+,e_k^+,e_j^-)-A(e_i^+,e_k^-,e_k^-,e_j^-)
    =-2\rho_2(e_i,e_j)=-2a_2\delta_{ij}$.
\medbreak\noindent This shows that $\rho e_i^\pm=2a_1e_i^\pm\pm 2a_2e_i^\mp$. Thus we may view $\rho$ as acting by
complex scalar multiplication by $\lambda:=2a_1+2a_2\sqrt{-1}$ on
$V_{\mathbb{C}}$. This implies that the underlying real operator is complex diagonalizable and
$\operatorname{Spec}(\rho)=\{\lambda,\bar\lambda\}$. Since the tensors
$A_i$ were extended to be complex multi-linear, we have
\begin{eqnarray*}
&&A(\rho v_1,v_2,v_3,v_4)=\operatorname{Re}\{A_1(\lambda v_1,v_2,v_3,v_4)+\sqrt{-1}A_2(\lambda
v_1,v_2,v_3,v_4)\}\\
&=&\operatorname{Re}\{\lambda A_1(v_1,v_2,v_3,v_4)+\sqrt{-1}\lambda A_2(v_1,v_2,v_3,v_4)\}\\
&=&\operatorname{Re}\{A_1(v_1,\lambda v_2,v_3,v_4)+\sqrt{-1}A_2(v_1,\lambda v_2,v_3,v_4)\}
=A(v_1,\rho v_2,v_3,v_4)\,.
\end{eqnarray*}
This establishes one equality of Lemma \ref{lem-1} (3); the other equalities
follow similarly and hence $\MM(\NN)$ is Jacobi--Ricci commuting as well.\hfill\qedbox

\medbreak We shall need the following technical result before establishing the second assertion of Theorem \ref{thm-6}.
Although well known, we include the proof for the sake of completeness and to establish notation:

\begin{lemma}\label{lem-8}
Let $J$ be a self-adjoint map of $(V,\langle\cdot,\cdot\rangle)$ so that
$J^2=-\operatorname{id}$. Then there exists an orthonormal basis $\{e_1^\pm,...,e_p^\pm\}$ for $V$
so that $Je_i^\pm=\pm e_i^\mp$.
\end{lemma}

\begin{proof} We assume that $p=1$ as the general result then follows by induction. Let $\{f^\pm\}$
be an orthonormal basis for $V$ where $f^+$ is spacelike and $f^-$ is timelike. As $J$ is trace-free and
self-adjoint,
$$J=\left(\begin{array}{rr}a&b\\-b&-a\end{array}\right)\,.$$
Since $J^2=-\operatorname{id}$, $b^2-a^2=1$. Let $e(\theta):=\cosh\theta f^++\sinh\theta
f^-$. Then 
\begin{eqnarray*}
&&\langle J e(\theta),e(\theta)\rangle\\
&=&\langle (a\cosh\theta+b\sinh\theta)f^++(-b\cosh\theta-a\sinh\theta)f^-,
\cosh\theta f^++\sinh\theta f^-\rangle\\
&=&a\cosh^2\theta+b\cosh\theta\sinh\theta+b\cosh\theta\sinh\theta+a\sinh^2\theta\\
&=&\textstyle\frac12\{(a+b)e^{2\theta}+(a-b)e^{-2\theta}\}\,.
\end{eqnarray*}
Since $a^2-b^2=-1$, $a+b$ and $a-b$ have opposite signs. Thus for some value of $\theta$, we have $\langle
Je(\theta),e(\theta)\rangle=0$. Set $e_i^+=e(\theta)$ and $e_i^-=Je(\theta)$.
\end{proof}

\medbreak\noindent{\it Proof of Theorem \ref{thm-6} (2):} Let $\MM=(V,\langle\cdot,\cdot\rangle,A)$ be a simple
Jacobi--Ricci commuting model. Assume that the Ricci operator $\rho$ is complex diagonalizable and that
$\operatorname{Spec}(\rho)=\{2a_1\pm2a_2\sqrt{-1}\}$ for
$a_2>0$. By Lemma
\ref{lem-8}, there is an orthonormal basis $\{e_i^\pm\}$ for $V$ so $Je_i^\pm=\pm e_i^\mp$. Set
$$
\begin{array}{l}
A_1(e_i,e_j,e_k,e_l):=A(e_i^+,e_j^+,e_k^+,e_l^+),\\
A_2(e_i,e_j,e_k,e_l):=-A(e_i^-,e_j^+,e_k^+,e_l^+)\,.
\end{array}$$
We may then derive the relations of Equations
(\ref{eqn-4}) from Lemma \ref{lem-4} (2). We check that $A_1$ and $A_2$ are algebraic curvature tensors
by verifying that:
\medbreak\qquad
$A_1(e_i,e_j,e_k,e_l)=A(e_i^+,e_j^+,e_k^+,e_l^+)=-A(e_j^+,e_i^+,e_k^+,e_l^+)$
\par\qquad$\quad=-A_1(e_j,e_i,e_k,e_l),$
\par\qquad$A_1(e_i,e_j,e_k,e_l)=A(e_i^+,e_j^+,e_k^+,e_l^+)=A(e_k^+,e_l^+,e_i^+,e_j^+)$
\par\qquad$\quad=A_1(e_k,e_l,e_i,e_j),$
\par\qquad$A_1(e_i,e_j,e_k,e_l)+A_1(e_j,e_k,e_i,e_l)+A_1(e_k,e_i,e_j,e_l)$
\par\qquad$\quad=A(e_i^+,e_j^+,e_k^+,e_l^+)+A(e_j^+,e_k^+,e_i^+,e_l^+)+A(e_k^+,e_i^+,e_j^+,e_l^+)=0$,
\par\qquad
$A_2(e_i,e_j,e_k,e_l)=-A(e_i^-,e_j^+,e_k^+,e_l^+)=A(e_j^+,e_i^-,e_k^+,e_l^+)$
\par\qquad$\qquad=-A(e_j^-,e_i^+,e_k^+,e_l^+)
=-A_2(e_j,e_i,e_k,e_l),$
\par\qquad$A_2(e_k,e_l,e_i,e_j)=-A(e_k^-,e_l^+,e_i^+,e_j^+)=-A(e_i^+,e_j^+,e_k^-,e_l^+)$
\par\qquad$\qquad=-A(e_i^-,e_j^+,e_k^+,e_l^+)=A_2(e_i,e_j,e_k,e_l),$
\par\qquad$A_2(e_i,e_j,e_k,e_l)+A_2(e_j,e_k,e_i,e_l)+A_2(e_k,e_i,e_j,e_l)$
\par\qquad$\qquad=-A(e_i^-,e_j^+,e_k^+,e_l^+)-A(e_j^-,e_k^+,e_i^+,e_l^+)-A(e_k^-,e_j^+,e_i^+,e_l^+)$
\par\qquad$\qquad=-A(e_i^-,e_j^+,e_k^+,e_l^+)-A(e_j^+,e_k^+,e_i^-,e_l^+)-A(e_k^+,e_j^+,e_i^-,e_l^+)=0$.
\medbreak\noindent
We verify $A_1$ and $A_2$ are Einstein by summing over $k$ to compute
\medbreak\qquad
$(\rho_1e_i,e_l)=A(e_i^+,e_j^+,e_j^+,e_l^+)$
\par\qquad\quad$=\frac12A(e_i^+,e_j^+,e_j^+,e_l^+)
-\frac12A(e_i^+,e_j^-,e_j^-e_l^+)=\frac12\langle\rho e_i^+,e_l^+\rangle=a_1\delta_{il}$.
\par\qquad
$(\rho_2e_i,e_l)=-A(e_i^-,e_j^+,e_j^+,e_l^+)$
\par\qquad\quad$=-\frac12A(e_i^-,e_j^+,e_j^+,e_l^+)
+\frac12A(e_i^-,e_j^-,e_j^-e_l^+)=-\frac12\langle\rho e_i^-,e_l^+\rangle=a_2\delta_{il}$.
\medbreak\noindent The desired result now follows.
\endproof

\medbreak Replacing $e_i^-$ by $-e_i^-$ in Equation
(\ref{eqn-4}) yields an isomorphism between the models $\NN(V,(\cdot,\cdot),A_1,A_2)$ and
$\NN(V,(\cdot,\cdot),A_1,-A_2)$; it is for this reason that we may always assume the Einstein constant of $A_2$ is
positive. This reflects that complex conjugation defines a field isomorphism of $\mathbb{C}$ taking
$\lambda\rightarrow\bar\lambda$ or, equivalently, by replacing $a_2$ by $-a_2$ in the construction. More
important, however, is the fact that the splitting
$V=V_+\oplus V_-$ where $V_\pm:=\operatorname{Span}\{e_i^\pm\}$ which is crucial to our discussion is highly
non-unique. Let $\NN=(V_0,g,A_1,A_2)$ and let $\tilde\NN=(\tilde V_0,\tilde g,\tilde A_1,\tilde A_2)$.
Let $\MM=\MM(\NN)$ and $\tilde\MM=\MM(\tilde\NN)$. Let $J$ and $\tilde J$ be the
associated complex structures on $V$ and on $\tilde V$, respectively.  We then have maximal spacelike subspaces
$V_+:=V_0$ and
$\tilde V_+:=\tilde V_0$ of
$V$ and
$\tilde V$, respectively, so that for all $x$, $y$, $z$, $w$ in $V_0$ and for all $\tilde x$, $\tilde y$, $\tilde
z$, $\tilde w$ in $\tilde V_0$,
$$
\begin{array}{ll}
V_+\perp JV_+,&\tilde V_+\perp\tilde J\tilde V_+,\\
 A_1(v,w,x,y)=A(v,w,x,y),
&\tilde A_1(\tilde v,\tilde w,\tilde x,\tilde y)=\tilde A(\tilde v,\tilde w,\tilde x,\tilde y),\\
A_2(v,w,x,y)=A(Jv,w,x,y),&
\tilde A_2(\tilde v,\tilde w,\tilde x,\tilde y)=-\tilde A(\tilde J\tilde v,\tilde w,\tilde x,\tilde y),\\
g(x,y)=\langle x,y\rangle,&
\tilde g(\tilde x,\tilde  y)=\langle\tilde x,\tilde  y\rangle\,.
\end{array}$$

Suppose that $\Theta$ is an isomorphism from $\MM$ to $\tilde\MM$. We may then identify $V=\tilde V$ and
$J=\tilde J$. The decomposition $V=V_+\oplus JV_+$ defines orthogonal projections $\pi_\pm$. Since $\tilde V_+$
is spacelike, $\pi_+$ defines an isomorphism $\theta$ from $\tilde V_+$ to $V_+$. Let
$$T=-J\circ\pi_-\circ\theta^{-1}:V_+\rightarrow\tilde V_+\rightarrow V_-\rightarrow V_+\,.$$
We may then represent any element of $\tilde V_+$ in the form $v+JTv$ for $v\in V_+$.
\begin{lemma}\label{lem-9}
Adopt the notation established above:
\begin{enumerate}
\item $\tilde V_+\perp J\tilde V_+$ if and only if $T$ is skew-adjoint.
\item The induced metric on $\tilde V_+$ is positive definite if and only if $|T|<1$.
\end{enumerate}\end{lemma}

\begin{proof} We have that $J$ is self-adjoint and that
$J^2=-\operatorname{id}$. Consequently, we have the following implications which establish Assertion (1):
\medbreak\qquad\qquad\phantom{.....}$\tilde V_+\perp J\tilde V_+$.
\smallbreak\qquad\qquad$\Leftrightarrow$ $\langle v+JTv,Jw+JJTw\rangle=0$ for all $v,w\in V_+$.
\smallbreak\qquad\qquad$\Leftrightarrow$ $-\langle v,Tw\rangle-\langle Tv,w\rangle=0$ for all $v,w\in V_+$.
\smallbreak\qquad\qquad$\Leftrightarrow$ $T$ is skew-adjoint.
\medbreak\noindent We argue similarly to prove Assertion (2):
\medbreak\qquad\qquad\phantom{......}$\langle v+JTv,v+JTv\rangle>0$ for all $0\ne v\in V_+$.
\smallbreak\qquad\qquad$\Leftrightarrow$ $g(v,v)-g(Tv,Tv)>0$ for all $0\ne v\in V_+$.
\smallbreak\qquad\qquad$\Leftrightarrow$ $|Tv|^2<|v|^2$ for all $0\ne v\in V_+$.
\smallbreak\qquad\qquad$\Leftrightarrow$ $|T|<1$.
\end{proof}

\medbreak\noindent{\it Proof of Theorem \ref{thm-6} (3).} Suppose
$\Theta:\MM(V_0,A_1,A_2)\rightarrow\MM(\tilde V_0,\tilde A_1,\tilde A_2)$ is an isomorphism. We use $\Theta$ to
identify $V$ with $\tilde V$ and to parametrize $\tilde V_+$ in the form $\{v+JTV\}$ where $T$ is a skew-adjoint
linear map of $V_0$ with $|T|<1$. We then the following identities for all $v$, $w$, $x$, and $y$:
\begin{equation}\label{eqn-5}
\begin{array}{l}
\tilde g(v,w)=\langle v+JTv,w+JTw\rangle,\\
\tilde A_1(v,w,x,y)=A(v+JTv,w+JTw,x+JTx,y+JTy),\\
\tilde A_2(v,w,x,y)=A(J(v+JTv),w+JTw,x+JTx,y+JTy)\,.
\end{array}\end{equation}
Lemma \ref{lem-1} (2) and Equation (\ref{eqn-5}) imply that Equation
(\ref{eqn-2}) holds. This establishes one implication of Theorem \ref{thm-6} (3). As the arguments are reversible, the
converse implication holds as well.\hfill\qedbox

\section*{Acknowledgments} Research of P. Gilkey partially supported by the
Max Planck Institute in the Mathematical Sciences (Leipzig) and by Project MTM2006-01432 (Spain). 
Research of S. Nik\v cevi\'c partially supported by Project 144032 (Srbija).

\end{document}